\input amstex
\documentstyle{amsppt}
\magnification=1200
\vsize19.5cm
\hsize13.5cm
\TagsOnRight
\pageno=1
\baselineskip=15.0pt
\parskip=3pt

\def\p{\partial}
\def\noo{\noindent}
\def\eps{\varepsilon}
\def\lam{\lambda}
\def\Om{\Omega}

\def\pom{{\p \Om}}

\def\R{\bold R}

\def\th{\theta}
\def\wtt{\tilde}
\def\back{\backslash}
\def\Ga{\Gamma}

\def\dist{\text{dist}}
\def\det{\text{det}}

\def\ol{\overline}

\def\D{\nabla}
\def\phi{\varphi}

\def\M{\Cal M}

\def\w{\wtt w}

\def\kk{k}

\nologo
 \NoRunningHeads

\topmatter

\title{On Harnack inequalities and singularities of\\
    admissible metrics in the Yamabe problem}\endtitle

\author{Neil S. Trudinger\ \ \ Xu-Jia Wang }\endauthor

\affil{ The Australian National University }\endaffil

\address
 Centre for Mathematics and its Applications, Australian National
    University, \newline Canberra ACT 0200, Australia
 \endaddress
\email neil.trudinger\@maths.anu.edu.au \ \ \
wang\@maths.anu.edu.au
\endemail

\thanks{
This work was supported by the Australian Research Council.
}\endthanks

\abstract { In this paper we study the local behaviour of
admissible metrics in the k-Yamabe problem on compact Riemannian
manifolds $(M, g_0)$ of dimension $n\ge 3$. For $n/2 <k<n$, we
prove a sharp Harnack inequality for admissible metrics when
$(M,g_0)$ is not conformally equivalent to the unit sphere $S^n$
and that the set of all such metrics is compact. When $(M,g_0)$ is
the unit sphere we prove there is a unique admissible metric with
singularity. As a consequence we prove an existence theorem for
equations of Yamabe type, thereby recovering a recent result of
Gursky and Viaclovski on the solvability of the $k$-Yamabe problem
for $k>n/2$. }\endabstract


\endtopmatter


\document

\baselineskip=13.2pt
\parskip=3pt

\centerline {\bf 1. Introduction}

\vskip10pt

Let $(\M, g_0)$ be a compact Riemannian manifold of dimension
$n\ge 3$ and $[g_0]$ the set of metrics conformal to $g_0$. For
$g\in [g_0]$ we denote by
$$A_{g}=\frac 1{n-2}(Ric_{g}-\frac {R_{g}} {2(n-1)}g)\tag 1.1$$
the Schouten tensor and by  $\lam(A_g)=(\lam_1, \cdots, \lam_n)$
the eigenvalues of $A_g$ with respect to $g$ (so one can also
write $\lam=\lam(g^{-1}A_g)$), where $Ric$ and $R$ are
respectively the Ricci tensor and the scalar curvature. We also
denote as usual
$$\sigma_k(\lam)
  =\sum_{i_1<\cdots<i_k}\lam_{i_1}\cdots\lam_{i_k}\tag 1.2$$
the $k$-th elementary symmetric polynomial and
$$\Ga_k=\{\lam\in \R^n\ |\
        \sigma_j(\lam)>0\ \text{for}\ j=1, \cdots, k\} \tag 1.3$$
the corresponding open, convex cone in $\R^n$. Denote
$$[g_0]_k=\{g\in [g_0]\ |\ \lam(A_g)\in\Ga_k\}. \tag 1.4$$
We call a metric in $[g_0]_k$ {\it $k$-admissible}. In this paper
we prove three main theorems pertaining to the cases $k>\frac n2$.

\proclaim{Theorem A} If $(\M, g_0)$ is not conformally equivalent
to the unit sphere $S^n$ and $\frac n2<k\le n$, then $[g_0]_k$ is
compact in $C^0(\M)$ and satisfies the following Harnack
inequality, namely for any $g=\chi g_0\in [g_0]_k$,
$$\max_{x, y\in\M}\frac{\chi(x)}{\chi(y)}
            \le exp(C|x-y|^{2-\frac nk}) \tag 1.5$$
for some fixed constant $C$ depending only on $(\M, g_0)$, where
$|x-y|$ denotes the geodesic distance in the metric $g_0$ between
$x$ and $y$.
\endproclaim

When the manifold $(\M, g_0)$ is the unit sphere, the compactness
is no longer true. In this case $(\M, g_0)$ is conformally
equivalent to the Euclidean space $\R^n$ so that without loss of
generality, it suffices to study conformal metrics on $\R^n$. For
our investigation we will allow singular metrics. Accordingly we
call a metric $g=\chi g_0$  $k$-admissible if $\chi:\ \M\to
(-\infty, \infty]$, $\chi$ is lower semi-continuous, $\not\equiv
\infty$ and there exists a sequence of $k$-admissible metrics
$g_m=\chi_m g_0$, $\chi_m\in C^2(\M)$, such that $\chi_m\to\chi$
almost everywhere in $\M$. If $g$ is $k$-admissible, then the
function $v=\chi^{(n-2)/4}$ is subharmonic with respect to the
operator
$$\square := -\Delta_g+\frac {n-2}{4(n-1)}R_g\tag 1.6$$
and hence by the weak Harnack inequality [GT], the set
$\{\chi=\infty\}$ has measure zero. Our next result classifies the
possible singularities of $k$-admissible metrics on $\R^n$.

\proclaim{Theorem B} Let $g$ be $k$-admissible on $\R^n$ with
$\frac n2<k\le n$. Then either
$$g(x)=\frac{C}{|x-x_0|^4} g_0(x)\tag 1.7$$
for some point $x_0\in\R^n$ and positive constant $C$, or the
conformal factor $\chi$ is H\"older continuous with exponent
$\alpha = 2-\frac nk$, where $g_0$ is the standard metric on
$\R^n$.
\endproclaim

\noo{\bf Remark}. Theorems A and B also hold if the condition
$g\in [g_0]_k$ (namely $\lam(A_g)\in\Ga_k$) is replaced by
$\lam(A_g)\in\Sigma_\delta$ for $\delta<\frac {1}{n-2}$, where the
cone
$$\Sigma_\delta=\{\lam\in\R^n\ |\ \lam_i >
      -\delta\sum_{j=1}^n \lam_j\ \ \forall\ \ 1\le i\le n\} \tag 1.8$$
was introduced in [GV2]. If $\lam\in\Ga_k$, then
$\lam\in\Sigma_\delta$ with $\delta=\frac {n-k}{n(k-1)}$ [TW2].

Theorems A and B have various interesting consequences. As an
application of Theorem A, we study the problem of prescribing the
$k$-curvature, that is the existence of a conformal metric $g\in
[g_0]$ such that
$$\sigma_k(\lam(A_g))=f, \tag 1.9$$
where $f$ is a given positive smooth function on $\M$. Write
$g=v^{4/(n-2)}g_0$. Then equation (1.9) is equivalent to the {\it
conformal $k$-Hessian equation}
$$\sigma_k(\lam(V))=\phi(x, v), \tag 1.10$$
where
$$V=-\D^2 v+\frac {n}{n-2}\frac {\D v\otimes\D v}{v}
    -\frac{1}{n-2}\frac{|\D v|^2}{v}g_0
    +\frac {n-2}{2} vA_{g_0},  \tag 1.11$$
$\lam(V)$ denotes the eigenvalues of the matrix $V$, and
$\phi=fv^{k\frac {n+2}{n-2}}$. When $f\equiv 1$, (1.9) is the {\it
$k$-Yamabe problem}, which has been studied by many authors, see
[A1,S, T] for $k=1$ and [CGY2, GeW, GW2, LL1, STW, GV1] for $k\ge
2$.

When $k\ge 2$, equation (1.10) is a fully nonlinear partial
differential equation, which is elliptic if the eigenvalues
$\lam(A_{g})\in \Ga_k$. Therefore to study problem (1.9), we
always assume $[g_0]_k\ne \emptyset$. Under this assumption, the
$\kk$-Yamabe problem has been solved in [STW] if $2\le k\le \frac
n2$ and (1.9) is variational. Equation (1.9) is automatically
variational when $k=2$, but when $k\ge 3$, it is variational when
the manifold is locally conformally flat or satisfies some other
conditions [STW].  When $\frac n2 <k \le n$, the existence of
solutions to (1.9) was proved in [GV1] for any smooth, positive
functions $f$;  see also [CGY2] for the solvability when $k=2$ and
$n=4$, and [GW2, LL1] when the manifold is locally conformally
flat. As a consequence of Theorem A, we have the following
stronger result.

\proclaim{Theorem C} Let $(\M, g_0)$ be a compact $n$-manifold not
conformally equivalent to the unit sphere $S^n$. Suppose $\frac
n2<k\le n$ and $[g_0]_k\ne \emptyset$. Then for any smooth,
positive function $f$ and any constant $p\ne k $, there exists a
positive solution to the equation
$$\sigma_k(\lam(V))=f(x)v^p. \tag 1.12$$
The solution is unique if $p<k$. When $p=k$, then there exists a
unique constant $\th>0$ such that
$$\sigma_k(\lam(V))=\th f(x)v^k \tag 1.13$$
has a solution, which is unique up to a constant multiplication.
\endproclaim

We may call the constant $\th$ in (1.13) (with $f\equiv 1$) the
{\it eigenvalue} of the conformal $k$-Hessian operator in (1.10).
As a special case of Theorem C, letting $p=k\frac {n+2}{n-2}$, we
obtain the existence of solutions to the $k$-Yamabe problem (1.9)
for $\frac n2<k\le n$, which was first proved in [GV1]. We also
include some extensions of Theorem C at the end of Section 4.

As in [STW] we will use conformal transforms of different forms,
$$g=\chi g_0=v^{\frac {4}{n-2}}g_0=u^{-2}g_0=e^{-2w}g_0\tag 1.14$$
so that
$$u=v^{-2/(n-2)}=e^w. \tag 1.15$$
We say $u$, $v$, or $w$ is {\it conformally $k$-admissible}, or
simply $k$-admissible if no confusion arises, if the metric $g$ is
$k$-admissible. In the smooth case, from the matrix $V$ in (1.11),
we see that $u, w$ are $k$-admissible if the eigenvalues of the
matrices
$$\align
U & =\{u_{ij}-\frac {|Du|^2}{2u}g_0 +uA_{g_0}\},\tag 1.16\\
W & =\{w_{ij}+w_iw_j -\frac 12 |Dw|^2 g_0 +A_{g_0}\}\tag 1.17\\
\endalign $$
lie in $\ol\Ga_k$, the closure of $\Ga_k$. Note that if $g$ is the
metric given by (1.7), then
$$v=\frac C{|x-x_0|^{n-2}}\tag 1.18$$
is the fundamental solution of the Laplace operator.

The conformal $k$-Hessian equation is closely related to the
$k$-Hessian equation
$$\sigma_k(\lam(D^2 u))=\phi\ \ \ \text{in}\ \ \Om,\tag 1.19$$
where $\Om \subset\R^n$ is a bounded domain. For the $k$-Hessian
equation (1.19), it is proved in [TW2] that when $\frac n2<k\le
n$, a $k$-admissible function (relative to  equation (1.19)) is
locally H\"older continuous with H\"older exponent $\alpha=2-\frac
nk$. The existence of solutions to (1.19) with right hand side
$\phi=f(x)|u|^p$ for some constant $p>0$ was studied in [CW] for
$k\le \frac n2$ and in [Ch, W] for $k=n$. By the H\"older
continuity one can extend the results in [Ch, W] to the cases
$\frac n2<k\le n$. The argument in [W] uses a degree theory, which
does not require a variational structure. We will employ the same
degree argument to prove our Theorem C.

We will first prove Theorem B for radially symmetric,
$k$-admissible functions defined on $\R^n$, then extend it to
general $k$-admissible functions by the comparison principle. The
proof of Theorem B also implies that if $w$ is a $k$-admissible
function on a manifold $\M$, then either $w$ is H\"older
continuous, or
$$w=-2\log |x-x_0|+C+o(1)\tag 1.20$$
for some point $x_0\in\M$. If the case (1.20) occurs, we show that
$w$ must be a smooth function. Hence by Bishop's volume growth
formula, it occurs only when the manifold is conformally
equivalent to the unit sphere, because when $\frac n2<k\le n$,
$\M$ equipped with the metric $g=e^{-2w}g_0$ is a complete
manifold with nonnegative Ricci curvature. Theorem C follows from
Theorem A and a degree argument.

The above theorems extend to more general symmetric curvature
functions. For example the $k^{th}$ elementary symmetric
polynomial $\sigma_k$ in (1.9) can be replaced by the quotient
$\sigma_k/\sigma_l$, where $k>l\ge 1$ and $n\ge k>\frac n2$. In a
subsequent paper we will extend these results to more general
symmetric curvature functions, as well as to the case $k=\frac n2$
in Theorem C.

\vskip30pt

\centerline{\bf 2. Proof of Theorem B}

\vskip10pt

\noo{\bf 2.1. Radial functions}. The proof of Theorem B can be
included in that of Theorem A. However we provide a separate proof
here. We first consider radially symmetric functions.  Let $w$ be
a radially symmetric, $k$-admissible function on
$\R^n\backslash\{0\}$. For any given point $x\ne 0$, by a rotation
of axes we assume $x=(0, \cdots, 0, r)$. Regard $w$ as a function
of $r=|x|$, $r\in (0, \infty)$. Then the matrix $W$ in (1.17) is
diagonal,
$$  W =\text{diag}(
  \frac 1r w'-\frac 12{w'}^2, \cdots,
  \frac 1r w'-\frac 12{ w'}^2,
   w''+\frac 12 {w'}^2) .$$
Denote $a=w''+\frac 12 {w'}^2$ and $b=\frac 1r w'-\frac 12{w'}^2$.
We have
$$\align
\sigma_k(\lam(W))
  & =b^kC_{n-1}^k+ab^{k-1} C_{n-1}^{k-1}\\
  & = C_{n-1}^{k-1}b^{k-1} (a+\frac{n-k}{k} b) .\tag 2.1\\
  \endalign $$
Since $\lam(W)\in\ol\Ga_k$ and $k>\frac n2$,
$$\align
 b & =\frac {w'}{r}-\frac 12 {w'}^2\ge 0,\tag 2.2\\
 a+\frac {n-k}{k} b &=(w''+\frac {w'}{r})
    -(1-\th)(\frac {w'}{r}-\frac 12 {w'}^2)\ge 0,\tag 2.3\\
    \endalign $$
where $\th=\frac {n-k}{k}<1$. It follows that
$$\align
& 0\le  w'  \le \frac 2r,\tag 2.4\\
& w''+\frac {w'}r \ge 0.\tag 2.5\\
\endalign $$
Note that (2.5) can also be written as $(rw')'\ge 0$. Therefore we
have

\proclaim{Lemma 2.1} The function $rw'$ is nonnegative, monotone
increasing, and $rw'\le 2$.
\endproclaim

It follows that $w$ must be locally uniformly bounded from above.
Next we prove

\proclaim{Lemma 2.2} The function $w$ is either H\"older
continuous in $\R^n$ with exponent $\alpha=2-\frac nk$, or
$$w(r)=2\log r+C\tag 2.6$$
for some constant $C$.
\endproclaim

\noo{\it Proof}. First we consider the case $k=n$. In this case
$a=w''+\frac 12{w'}^2\ge 0$, namely, $\frac{w''}{{w'}^2}+\frac
12\ge 0$. Hence
$$\int_0^r(\frac {-1}{w'}+\frac r2)'\ge 0.$$
If $w$ is not Lipschitz continuous, we have $w'(r)\to\infty$ as
$r\to 0$. Hence
$$\frac {-1}{w'}+\frac r2\ge 0. $$
It follows that $w' \ge\frac 2r$. Hence  by Lemma 2.1, $w'\equiv
2/r$ so that $w(r)=2\log r +C$.

In the cases $\frac n2<k<n$, if $rw'\not\equiv 2$, then by Lemma
2.1,  $\lim_{r\to 0} rw'=c_0<2$. For any $c_1\in (c_0, 2)$,
$$w''+\frac {w'}r \ge (1-\th)\frac {w'}{r} (1-\frac 12rw')\ge
  (1-\th)(1-\frac{c_1} 2)\frac {w'}{r} \tag 2.7$$
if $r$ is sufficiently small. Hence
$$\frac {w''}{w'}+ \frac\sigma r\ge 0,$$
where $\sigma=1-(1-\th)(1-\frac {c_1} 2)<1$.  We obtain
$$\log (w'r^{\sigma})\big|^{r_0}_r\ge 0.$$
Hence
$$w'\le \frac {C}{r^{\sigma}}. \tag 2.8$$
Hence $w$ is bounded and continuous.

To show that $w$ is H\"older continuous with H\"older exponent
$\alpha=2-\frac nk$, by Lemma 2.1 it suffices to prove it at
$r=0$. Note that
$$a+\th b=w''+\th\frac {w'}{r}+\frac{1-\th}{2}{w'}^2\ge 0.$$
Hence
$$\frac {w''}{w'}+\frac \th r\ge -\frac{1-\th}{2}w'.$$
Taking integration from $r$ to $r_0$, we obtain
$$\log (w'r^\th)\big|^{r_0}_r\ge C. $$
Hence
$$w'\le \frac {C}{r^{\th}}, \tag 2.9$$
so that $w$ is H\"older continuous with exponent $1-\th=2-\frac
nk$. $\square$

\noo{\it Remark 2.1}. The H\"older continuity also follows from
[TW2]. Let $u=e^w$ as in (1.15). Then from the matrix $U$ in
(1.16) we see that $u$ is $k$-admissible with respect to the
$k$-Hessian operator $\sigma_k(\lam(D^2 u))$. Hence $u$ is
H\"older continuous with exponent $\alpha=2-\frac nk$. It follows
that for any constant $c>0$, $w_c=\max(w, -c)$ is also H\"older
continuous with exponent $2-\frac nk$. In particular, if $w_m$
converges to $w$ a.e., then $w_m$ converges to $w$ uniformly in
$\{w>-c\}$ for any $c>0$.

\vskip10pt

\noo{\bf 2.2. Proof of Theorem B}. Let $w$ be a $k$-admissible
function. For any $h\in\R$, denote $\Om_h=\{w<h\}$. Since $w$ is
upper semi-continuous, $\Om_h$ is an open set. For any given point
$0$, we define a function $\wtt w$ of one variable $r$ by
$$\wtt w(r)=\inf \{h: \ \dist (0, \pom_h)>r\}.\tag  2.10$$
Let $x_h\in\pom_h$ such that $|x_h|=r_h:=\dist (0, \pom_h)$.
Assume that $\pom_h$ and $w$ are smooth at $x_h$.  Rotate the axes
such that $x_h=(0, \cdots, 0, r_h)$. Then the $x_n$-axis is the
outer normal of $\pom_h$ at $x_h$. Hence
$$\align
\wtt w(r_h) & =w(x_h),\\
\wtt w(r_h+t) &\ge w(x_h+te_n)\tag 2.11\\
\endalign $$
for $t$ near $0$, where $e_n=(0, \cdots, 0, 1)$. We obtain
$$\align
\w'(r_h)& =w_n(x_h)=|Dw|(x_h) ,\tag 2.12 \\
\w''(r_h) &\ge w_{nn}(x_h)  \\
\endalign $$
provided $\w$ is twice differentiable point at $r_h$.

Let $\kappa_1, \cdots, \kappa_{n-1}$ be the principal curvatures
of $\pom_h$ at $x_h$. Then
$$w_{ij}=|Dw|\kappa_i\delta_{ij}\ \ \ i, j\le n-1. \tag 2.13$$
By our choice of $x_h$, we have
$$\kappa_i\le \frac 1r,\tag 2.14$$
where $r=r_h$. Hence the matrix
$$(w_{ij})_{i, j=1}^{n-1}\le \frac 1r |Dw| I. \tag 2.15$$
At $x_h$, the matrix $W$ is given by
$$\align
W & =\{w_{ij}+w_iw_j -\frac 12 |Dw|^2 I\}\\
 &=\pmatrix
  w_{11}-\frac 12|Dw|^2, &0, & \cdots, & w_{1n}\\
  0, & w_{22}-\frac 12 |Dw|^2, & \cdots, & w_{2n}\\
  \cdot &\cdot &\cdot & \cdot\\
  \cdot &\cdot &\cdot & \cdot\\
  w_{1n}, & w_{2n}, &\cdots, & w_{nn}+\frac 12 |Dw|^2\\
  \endpmatrix. \endalign $$
Let
$$W' =\text{diag}(
  w_{11}-\frac 12|Dw|^2,\cdots, w_{22}-\frac 12 |Dw|^2,
  w_{nn}+\frac 12 |Dw|^2)\tag 2.16$$
be a diagonal matrix.  We claim that the eigenvalues
$\lam(W')\in\ol\Ga_k$. Indeed, recalling that $\sigma_k(\lam(W))$
is the sum of all principal $k\times k$ minors, we have
$$\sigma_k(\lam(W))=\sigma_k(\lam(W'))
  -\sum_{i<n} \sigma_{k-2}(\lam(W_{|in}))w_{;in}^2, \tag 2.17$$
where $w_{;ij}$ is the entry of the matrix $W$, and $W_{|ij}$
denotes the matrix obtained by cancelling the $i$th and $j$th rows
and columns of $W$. Since $\lam(W)\in\ol\Ga_k$, we have
$$\sigma_{k-2}(\lam(W_{|in}))
  =\frac {\p^2 \sigma_k(\lam(W))} {\p w_{;ii} \p w_{;nn}} >0 .\tag 2.18$$
Hence $\sigma_k(\lam(W'))\ge \sigma_k(\lam(W))\ge 0$. Similarly we
have $\sigma_j(\lam(W'))\ge \sigma_j(\lam(W))$ for $1\le j \le k$,
and so $\lam(W')\in\ol\Ga_k$.

From (2.15),
$$W'\le \text{diag} (
  \frac 1r \wtt w'-\frac 12 (\wtt w')^2, \cdots,
   \frac 1r \wtt w'-\frac 12(\wtt w')^2,
     \wtt w''+\frac 12 (\wtt w')^2) .\tag 2.19$$
Therefore as in \S 2.1, we see that $\w$ satisfies
$$\align
& \frac {\wtt w'}{r}-\frac 12(\wtt w')^2\ge 0\tag 2.20\\
& (\w''+\frac {\wtt w'}{r})
    -(1-\th)(\frac {\wtt w'}{r}-\frac 12(\wtt w')^2\ge 0 \tag 2.21\\
    \endalign $$
if $\w$ is twice differentiable at $r$.

To proceed further we need some remarks.

\noo{\it Remarks 2.2}.
\newline
(i) If the function $\w$ is not smooth, by (2.11) it satisfies
(2.20) and (2.21) in the viscosity sense. That is if $\phi$ is a
smooth function satisfying
$$\align
& \frac {\phi'}{r}-\frac 12 {\phi'}^2\ge 0,\\
& (\phi''+\frac {\phi'}{r})
    -(1-\th)(\frac {\phi'}{r}-\frac 12(\phi')^2= 0,\\
    \endalign $$
and $\w(r_0)=\phi(r_0)$, $\w'(r_0)=\phi'(r_0)$, then $\w(r)\ge
\phi(r)$ near $r_0$. If instead $\w(r_0)=\phi(r_0)$,
$\w(r_1)=\phi(r_1)$, then $\w(r)\le \phi(r)$ for $r\in (r_0,
r_1)$.
\newline
(ii) In the above we assumed that both $w$ and $\pom_h$ are smooth
at $x_h$. If $w$ is smooth but $\pom_h$ is not smooth at $x_h$, it
is easy to see that (2.15) still holds and so one also has (2.20)
and (2.21). If $w$ is not smooth, by definition it can be
approximated by smooth functions. Hence (2.20) and (2.21) always
hold.
\newline
(iii) Another way to verify (2.20) and (2.21) is to regard $\w$ as
a function of $x$, namely $\w(x)=\w(|x|)$. Then $\w-w$ attains a
local minimum at $x_h$. Hence $\w$ is $k$-admissible in the
viscosity sense, and so (2.20) and (2.21) hold.

From (2.20) and (2.21), we can prove Theorem B easily. First we
consider the case when $w$ is unbounded from below.

\proclaim{Lemma 2.3} Let $w$ be a $k$-admissible function which is
unbounded from below, then there exists a point $x_0\in\R^n$ and a
constant $C$ such that
$$w(x)\equiv -2\log |x-x_0|+C.\tag 2.22$$
\endproclaim

\noo{\it Proof}. If $w$ is unbounded from below, the singular set
$S=\bigcap_{\{c<0\} }\{w<c\}$ is not empty. Choose a point $0\in
S$.  By (2.20) and (2.21), and from the argument in \S 2.1, we
must have $\w(r)=2\log r+C$ for some constant $C$.

Let $\hat w=2\log |x|+C$. Then
$$\align
& \sigma_1(\lam(W_{\hat w}))=0,\\
& \sigma_1(\lam_1(W_w))\ge \sigma_k^{1/k}(\lam(W_w))\ge 0,\\
\endalign $$
where $W_{\hat w}$ is the matrix corresponding to $\hat w$, given
in (1.17). By the relation (1.15), $\sigma_1(\lam(W))$ is indeed
the Laplace operator. Since $\w=2\log r+C$, we see that $w-\hat w$
attains its local maximum at some interior point. By the maximum
principle for the Laplace equation, we conclude that $w\equiv \hat
w$. $\square$

Next we consider the case when $w$ is bounded from below.

\proclaim{Lemma 2.4} Let $w$ be a $k$-admissible function $w$.
Suppose $w$ is bounded from below. Then $w$ is H\"older continuous
with exponent $\alpha = 2-\frac nk$.
\endproclaim

\noo{\it Proof}. For any given point $x_0$, we may take $x_0$ as
the origin and define $\w$ as (2.10). Then to prove that $w$ is
H\"older continuous at $x_0$ with exponent $\alpha = 2-\frac nk$,
it suffices to show that $\w$ is H\"older continuous with exponent
$\alpha$. But by (2.20), (2.21), the H\"older continuity of $\w$
readily follows from the argument in \S 2.1, see (2.9). $\square$

\vskip5pt

The H\"older continuity also follows from Remark 2.1 above.

Note that the function $w=2\log |x|$ is $k$-admissible. By
truncating at $w=-K$ (for large $K$) and capping off, we see that
the set of H\"older continuous $k$-admissible functions is not
compact.

\vskip10pt

\noo{\bf 2.3. Applications}. First we remark that, by the above
proof, Theorem B also holds for $k$-admissible functions defined
on a domain. Here we restate the theorem for the function
$v=e^{-\frac{n-2}2 w}$. Note that by Lemma 2.1, a (non-smooth)
$k$-admissible function $v$ must be locally strictly positive when
$k>\frac n2$.

\proclaim{Theorem B$'$} Let $\Om$ be a domain in $\R^n$. Let $v$
be a $k$-admissible function in $\Om$ with $\frac n2<k\le n$. If
$v$ is unbounded from above near some point $x_0\in\Om$, then
$$v(x)=C|x-x_0|^{2-n}.\tag 2.23 $$
Otherwise $v$ is locally H\"older continuous in $\Om$ with
exponent $\alpha = 2-\frac nk$.
\endproclaim

It was proved in [LL1] that if $v$ is a $k$-admissible function,
so is the function $v_\psi$ in $B_1(0)\back\{0\}$, where
$$v_\psi=|J_\psi|^{\frac {n-2}{2n}} v\cdot\psi\tag 2.24$$
$\psi(x)=\frac {x}{|x|^2}$, and $J_\psi$ is the Jacobian of the
mapping $\psi$. From Theorem B we have

\proclaim{Corollary 2.5} Let $v$ be a $k$-admissible function
defined in $\R^n\back B_1(0)$ with $\frac n2<k\le n$. Then either
$v\equiv constant$ or $|x|^{n-2}v(x)$ converges to a positive
constant as $x\to\infty$.
\endproclaim

\noo{\it Proof}. We cannot apply Theorem B$'$ directly, as the
function $v_\psi$ has a singular point at $0$. Denote $w=\frac
{-2}{n-2}\log v_\psi$. If $w(x)\to -\infty$ as $x\to 0$, the
argument in \S 2.2 implies that $w=2\log |x|+C$ and so $v\equiv
constant$. Otherwise it suffices to show that $w$ is continuous at
$0$.

Let $w(0)={\overline\lim}_{x\to 0}w(x)$ so that $w$ is upper
semi-continuous. If $a=:\underline{\lim}_{x\to 0} w(x)<w(0)$, for
simplicity let us assume that $a\le -1$ and $w(0)=0$. Let $x_m\to
0$ such that $w(x_m)=-1$. Define the function $\w=\w_{x_m}$ as in
(2.10), with center at $x_m$. We claim that when $m$ is
sufficiently large, the point $x_h$ in (2.11) at $h=0$ cannot be
the origin. Indeed, if $x_h=0$, by the H\"older continuity of $\w$
(in the range $-1<\w<0$) we see that $w(x)\le -\frac 12$ when
$|x-x_m|\le \delta |x_m|$ for some $\delta>0$ independent of $m$.
But note that $v_\psi=e^{-\frac{n-2}{2}w}$ is supharmonic.
Applying the mean value theorem to $e^{-\frac{n-2}{2}w}$ we
conclude that ${\overline\lim}_{x\to 0}w(x)>0$. This is a
contradiction.

It follows by the argument in \S 2.2 that $\w=\w_{x_m}$ is
uniformly H\"older continuous. Hence if $w(0)=0$ and $w(x_m)\le
-1$, we have $|x_m|\ge c_0>0$ for some $c_0$ independent of $m$.
This is again a contradiction. Hence $w$ is continuous at $0$, and
so $|x|^{n-2}v(x)$ converges to a positive constant as
$x\to\infty$. $\square$

By Theorem B$'$, we have either $v_\psi=2\log|x|+C$, or $v_\psi$
is H\"older continuous at $0$. Hence the results in Corollary 2.5
follows. Theorem B also implies the non-existence of solutions to
the Dirichlet problem in general. Let $\Om$ be a non-round,
bounded domain in $\R^n$ containing the origin.  Then if $k>\frac
n2$, there is no solution to the Dirichlet problem
$$\align
\sigma_k(\lam(V))& =f\ \ \ \text{in}\ \ \Om, \tag 2.25\\
v & =c\ \ \ \text{on}\ \ \pom\\
\endalign $$
in general, where  $c$ is any positive constant, and $f$ is a
positive smooth function. Indeed, let $\{f_m\}$ be a sequence of
smooth, positive functions which converges to zero locally
uniformly in $\Om\back\{0\}$ such that $\sup v_m\to\infty$, where
$v_m$ is the corresponding solution. Then $v_m$ must converge to
the function $v=C|x|^{2-n}$ by Theorem B. Hence $\Om$ must be a
ball.

For the existence of solutions to the Dirichlet problem, it was
proved in [G] that for any smooth, bounded domain with smooth
boundary data, if there exists a sub-solution, then there exists a
solution to the Dirichlet problem.

\vskip30pt

\newpage

\centerline{\bf 3. Proof of Theorem A}

\vskip10pt

\noo{\bf 3.1. H\"older continuity}. We start with a H\"older
continuity property of $k$-admissible functions.

\proclaim{Lemma 3.1} Let $(\M, g_0)$ be a compact manifold.
Suppose $g=u^{-2}g_0\in [g_0]_k$ and $k>\frac n2$. Then $u$ is
H\"older continuous with exponent $\alpha=2-\frac nk$,
$$\frac {u(x)-u(y)}{|x-y|^\alpha}\le C\int_\M u,\tag 3.1$$
where $C$ is independent of $u$.
\endproclaim

\noo{\it Proof}. By approximation it suffices to prove (3.1) for
smooth functions.  For any given point $0\in\M$, there exists a
conformal metric [A2,C,Gu], still denoted by $g_0$, such that in
the normal coordinates at $0$,
$$\det (g_0)_{ij}\equiv 1\ \ \ \text{near}\ \ 0. \tag 3.2$$
Let
$$u_0(x)=|x|^{2-\frac nk},\tag 3.3$$
where $|x|$ denotes the geodesic distance from $0$. Note that
under condition (3.2), the Laplacian $\Delta$ on $\M$ is equal to
the Euclidean Laplacian when applying to functions of $r=|x|$
alone [LP, SY]. Hence
$$\Delta_{g_0} u_0=\frac {n(k-1)(2k-n)}{k^2} r^{-\frac nk}. \tag 3.4$$
Denote by
$$P[u]=\min \lam_i+ \delta \sum_i\lam_i,
     \ \ \ (\delta=\frac{n-k}{n(k-1)})\tag 3.5$$
the Pucci minimal operator [GT], where $(\lam_1, \cdots, \lam_n)$
are the eigenvalues of the Hessian matrix $(\D_{ij} u_0)$.
Obviously we have
$$\min \lam_i\le\p_r^2 u_0
     =-\frac{(2k-n)(n-k)}{k^2} r^{-\frac nk}. \tag 3.6$$
Therefore $u_0$ satisfies
$$P[u_0]\le 0\ \ \text{in} \ \ B_{0, r}\back \{0\}. $$
where $B_{y,r}$ denotes the geodesic ball with center $y$ and
radius $r$.

On the other hand, since $\lam(U)\in\ol\Ga_k$, where $U$ is given
in (1.16), we have $\lam(u_{ij}+uA_{g_0})\in \ol\Ga_k\subset
\ol\Ga_1$. Namely $\Delta u+\text{tr}(A_{g_0})u\ge 0$. By the
Harnack inequality it follows
$$\sup u\le C\int_\M u .\tag 3.7$$
Therefore to prove (3.1) we may assume that $\int_\M u=1$ and $u$
is uniformly bounded.

Let $u_a=u+a |x|^2$. Then $\D^2 u_a> \D^2 u + aI$ near $0$, where
$I$ is the unit matrix.  Since $\lam(\D^2 u+uA_{g_0})\in
\ol\Ga_k$, we have $\lam(\D^2 u_a)\in\Ga_k$ when $a$ is suitably
large. Taking $l=1$ in the proof of Lemma 4.2 in [TW2], one has
$$\lam_i+\frac {n-k}{n(k-1)}\sum_i\lam_i\ge 0, \tag 3.8$$
namely $P[u_a]\ge 0$ near $0$. Hence by applying the comparison
principle to the functions $u_a$ and $u_0$ with respect to the
operator $P$, we conclude the H\"older continuity (3.1). $\square$

\noo{\it Remark}. The estimate (3.1) (with exponent
$\alpha<2-\frac nk$) also follows from gradient estimates from our
reduction to $p$-Laplacian subsolution in [TW2]. Since
$\lam(U)\in\Ga_k$, we have $\lam(D^2 u+uA_{g_0})\in\Ga_k$. By
(3.8) it follows that
$$\Delta_p u:=\D_i(|\D u|^{p-2} \D_iu)\ge -Cu|\D u|^{p-2}\tag 3.9$$
for $p-2=\frac {n(k-1)}{n-k}$ and some constant $C$. From our
argument in [TW2], we obtain $\int_\M |\D u|^q\le C$ for any
$q<nk/(n-k)$, whence by the Sobolev inequality, we infer (3.1) for
$\alpha<2-\frac nk$; (see also [GV2]).

By the relation $u=e^w$, we have the following

\proclaim{Corollary 3.2} Let $w$ be a $k$-admissible function.
Suppose $w\le 0$. Then for any $K>0$, there exists $C=C_K>0$,
independent of $w$,  such that when $w(y)>-K$,
$$\frac {w(x)-w(y)}{|x-y|^\alpha}\le C.\tag 3.10$$
\endproclaim

From (3.10), we see that if $w(x)\le -K-1$, then $|x-y|\ge
C_{K+1}^{1/\alpha}$. Also note that in Corollary 3.2, if we assume
that $w\le 0$ in $B_{y,r}$, then (3.10) holds for $x, y\in
B_{y,r/2}$ for some $C$ depending on $r$.

\vskip10pt

\noo{\bf 3.2. Singularity behaviour of $k$-admissible functions}.
Suppose $w$ is a $k$-admissible function. At any given point
$0\in\M$, we choose a conformal normal coordinate such that (3.2)
holds. In the conformal metric, the Ricci curvature vanishes at
$0$ [LP, SY]. Hence
$$|A_{g_0}|\le Cr\ \ \ \text{near}\ \  0. \tag 3.11$$
Define $\w$ as in (2.10). Then the argument thereafter is still
valid, except that (2.14) should be replaced by $\kappa_i\le \frac
1r +C$.  Hence from (2.19), we have
$$(\wtt b,\ \cdots,\ \wtt b,\ \wtt a) \in \ol\Ga_k, \tag 3.12$$
where
$$\align
 \wtt b& = (\frac 1r+C) \w'-\frac 12(\w')^2+Cr,\\
 \wtt a & = \w''+\frac 12 (\w')^2+Cr. \\
  \endalign$$
Hence similarly to (2.2) (2.3), we have $\wtt b  \ge 0$ and
$$ \wtt a+\frac {n-k}k\wtt b
 =[\w''+(\frac 1r+C)\w'+Cr]
   -(1-\th)[(\frac 1r+C)\w'-\frac 12(\w')^2+Cr]\ge 0.$$
It follows, similarly to (2.4) and (2.5),
$$\align
& \w'  \le \frac 2r+\frac {Cr}{\w'}+C,\tag 3.13\\
& \w''+(\frac 1r+C)\w'+Cr\ge 0.\tag 3.14\\
\endalign $$
From (3.13),
$$\w'  \le \frac 2r+C$$
for a different $C$. Therefore by (3.14), we obtain
$$(r\w')'+C\ge 0. $$
It follows that $r\w'+Cr$ is increasing. By the compactness of
$\M$, a $k$-admissible function $w$ must be bounded from above.

If $r\w'<2$ near $r=0$, then similarly to (2.7) (2.8), $\w$ is
bounded and H\"older continuous.

If $r\w'\to 2$ as $r\to 0$, then $r\w'+Cr\ge 2$, namely $\w'\ge
\frac 2r -C$.  Hence we obtain
$$\frac 2r+C\ge \w'\ge \frac 2r-C. \tag 3.15$$
We obtain
$$\w(r)=2\log r+C'+O(r).\tag 3.16$$
By subtracting a constant we assume that $C'=0$.

\proclaim{Lemma 3.3} If $\w$ satisfies (3.16), then near $0$,
$$w(x)=2\log |x|+o(1).\tag 3.17$$
\endproclaim

\noo{\it Proof}. We prove (3.17) by a blow-up argument. In a
normal coordinate system at $0$, let $y=c_mx$ and
$w_m(y)=w(x)+2\log c_m$, where $c_m$ is any sequence converging to
infinity. Let $\w_m$ be the corresponding function of $w_m$. Then
by (3.16),
$$\w_m(r)=2\log r+O(c_m^{-1}).\tag 3.18$$
Hence $\w_m\to 2\log r$.

For any fixed $r_0>0$ small, let $w_m(y_m)=\w_m(r_0)$
($|y_m|=r_0$). We may assume that $y_m\to y_0$. By the H\"older
continuity (Corollary 3.2), we may also assume that in a
neighborhood of $y_0$, $w_m$ converges uniformly to $w_\infty$.
Then $w_\infty$ is a $k$-admissible function defined on $\R^n$.
The comparison principle argument of Lemma 2.3 implies that
$w_\infty\equiv 2\log r$ in a neighborhood of $y_0$. The H\"older
continuity in Corollary 3.2 implies that if $w_\infty=2\log r$ at
some point, $w_\infty$ is well-defined nearby. The comparison
principle then implies that $w_\infty\equiv 2\log r$ near the
point. Hence $w_\infty\equiv 2\log r$ in $\R^n\back\{0\}$ and
(3.17) is proved. $\square$

From the proof of Lemma 3.3, we see that $w$ has only isolated
singularities. For if there is a sequence of singular points
$x_m\in\M$ which converges to a point $0$, we may choose
$c_m=|x_m|$ in the above argument. Then the limit function
$w_\infty$ has at least two singular points $0$ and $x^*=\lim
x_m/|x_m|$. To see that $x^*$ is a singular point of the limit
function $w_\infty$, we notice that the constant $C'$ is uniformly
bounded from above if $w$ is negative in a neighbourhood of $0$,
which in turn implies that $\lim_{x\to x^*} w_\infty
(x^*)=-\infty$. But the above argument shows that $w_\infty=2\log
r$. This is a contradiction. Next we show that $w$ has at most one
singular point.

\proclaim{Lemma 3.4} Let $w$ be a $k$-admissible function. Then
the singularity set
$$S_w=\bigcap_{h<0} \{x\in\M\ |\ w(x)<h\}\tag 3.19$$
contains at most one point.
\endproclaim

\noo{\it Proof}. If $S_w$ is not empty, it consists of finitely
many isolated points.  Let $g=e^{-2w}g_0$. By Lemma 3.3, $(\M\back
S_w, g)$ is a complete manifold with finitely many ends. Now
fixing a point $y\not\in S_w$, we consider the ratio
$$Q(r)=\frac {Vol(B_{y, r})}{r^n},\tag 3.20$$
where $B_{y,r}=B_{y, r}[g]$ is the geodesic ball of $(\M, g)$. By
definition,  there is a sequence of smooth $k$-admissible
functions $w_m$ which converges to $w$ locally uniformly. It is
easy to verify that for any fixed $y$ and $r$,  $Vol(B_{y,
r}[g_m])\to Vol(B_{y, r}[g])$ as $m\to\infty$, where
$g_m=e^{-2w_m}g_0$. From [GVW], the Ricci curvature of $(\M, g_m)$
is positive. Hence by the Bishop Theorem, the ratio
$Q_m(r)=Vol(B_{y, r}[g_m])/r^n$ is decreasing for all $m$. Sending
$m\to\infty$, we see that $Q$ is non-increasing in $r$. Hence
$$Q(0)\le \lim_{r\to 0} Q(r)\le \frac 1n \omega_n,\tag 3.21$$
where $\omega_n$ is the area of the unit sphere $S^{n-1}$.

On the other hand, denote $A_{r_1, r_2}=B_{0, r_2}[g_0]-B_{0,
r_1}[g_0]$, where $r_2>r_1>0$ are sufficiently small. We identify
$A_{r_1, r_2}$ with the Euclidean annulus $A^e_{r_1,
r_2}=\{x\in\R^n\ |\ r_1<|x|<r_2\}$ by the exponential map. By the
asymptotic (3.17), the volume of $A_{r_1, r_2}$ in the metric
$g=e^{-2w}g_0$ is a lower order perturbation of that in the metric
$g'=e^{-2w'}g_0$, where $w'=2\log |x|$. But in our normal
coordinates at $0$, by (3.2) the volume of $A_{r_1, r_2}$ in $g'$
is the same as that of $A^e_{r_1, r_2}$ with the metric
$g'_e=e^{-2w'}g_e$, where $g_e$ is the standard Euclidean metric.
Hence $\text{Vol}_{g'} A_{r_1, r_2}=\frac 1n
\omega_n(r_1^{-n}-r_2^{-n})$. Therefore as $r\to \infty$, each end
of the metric $g$ will contribute to the ratio $Q(r)$ a factor
$\frac 1n\omega_n$. Therefore we obtain
$$\lim_{r\to\infty} Q(r)= \frac mn \omega_n, \tag 3.22$$
where $m$ is the number of singular points of $w$. From (3.21) and
(3.22) we see that if $S_w$ is not empty, then $m$ must be equal
to $1$, namely $S_w$ is a single point. $\square$

\vskip10pt

\noo{\bf 3.3. Smoothness of $k$-admissible functions}. In this
subsection we prove the following smoothness result.

\proclaim{Lemma 3.5} Let $w$ be a $k$-admissible function $w$ with
a singular point $0$. Then $w$ is $C^\infty$ smooth away from $0$.
\endproclaim

\noo{\it Proof}. First we prove
$$\sigma_k(\lam(A_g))\equiv 0 \ \ \text{in}\ \ \M\back\{0\},\tag 3.23$$
where $g=e^{-2w}g_0$.  It suffices to prove that for any given
point $x_0\ne 0$ and a sufficiently small $r>0$ ($r<\frac 14|x|)$,
(3.23) holds in $B_{x_0, r}=B_{x_0, r}[g_0]$.

By definition, there exists a sequence of smooth $k$-admissible
functions which converges to $w$ in $B_{x_0, 2r}$ uniformly. Let
$\phi_m$ be the solution of the Dirichlet problem [G]
$$\align
& \sigma_k(\lam(A_{g_{\phi_m}}))=\eps_m
                 \ \ \ \text{in}\ \ B_{x_0, r},\tag 3.24\\
& \phi_m=w_m\ \ \ \ \text{on}\ \ \p B_{x_0, r},\\
\endalign $$
where $g_{\phi_m}=e^{-2\phi_m}g_0$, and $\eps_m$ is a small
positive constant such that $\sigma_k(\lam(A_{g_{w_m}}))>\eps_m$
($g_{w_m}=e^{-2w_m}g_0$). By the comparison principle we have
$\phi_m\ge w_m$ in $B_{x_0, r}$. Let $\hat w_m=w_m$ in $\M-B_{x_0,
r}$ and $\hat w_m=\phi_m$ in $B_{x_0, r}$. Then $\hat w_m$ is
$k$-admissible (see Corollary 3.8 below). Let $\hat w =
\lim_{m\to\infty} \hat w_m$. Then $\hat w$ is a $k$-admissible
function with singularity point $0$. Define the metric $\hat
g=e^{-2\hat w}g_0$ and the ratio $\hat Q(r)=\frac {Vol(B_{y,
r}[\hat g])}{r^n}$. Then from the proof of Lemma 3.4, we also have
$\hat Q\equiv \frac 1n \omega_n$.

To prove (3.23) it suffices to show that $\hat w\equiv w$. Noting
that $\hat w=w$ in $\M-B_{x_0, r}$ and $\hat w\ge w$ in $B_{x_0,
r}$, we have $B_{y, r}[\hat g]\supset B_{y, r}[g]$ for any $r>0$
and $y\ne 0$. If there exists a point $y\in B_{x_0, r}$ such that
$\hat w> w$ at $y$, then there exists a positive constant
$\delta>0$ such that for any $r>1$,
$$B_{y, r}[\hat g]\supset B_{y, r+\delta}[g] . $$
But this is impossible as both the ratios $Q(r)$ and $\hat Q(r)$
are constant.

By the interior second order derivative estimate in [GW1, STW], we
see that $w$ is $C^{1,1}$ smooth. Next we prove that $w$ is
$C^\infty$ smooth away from $0$. By the regularity of linear
elliptic equations [GT], it suffices to prove that
$v=w^{-\frac{n-2}{2} w}\in C^{1,1}$ is a strong solution to the
uniformly elliptic equation
$$-\Delta_{g_0} v+\frac {n-2}{4(n-1)}R_{g_0}v=0
    \ \ \text{in}\ \ \M\back\{0\},\tag 3.25$$
where $R$ is the scalar curvature of $(\M, g_0)$. Namely the
scalar curvature of $g=e^{-2w}g_0$ vanishes identically.

Equation (3.25) is not hard to prove, see \S 7.6 in [GV1]. Here we
provide a proof for completeness. Since $w\in C^{1,1}$, it is
twice differentiable almost everywhere.  Suppose at a point $0$,
$w$ is twice differentiable and the scalar curvature $R>0$. Then
with respect to normal coordinates of $g$ at $0$, we have the
expansion
$$\det g_{ij}=1-\frac 13 R_{ij}x_ix_j+o(|x|^2),\tag 3.26 $$
see (5.2) in [LP]. Hence
$$\align
 \text{Vol}(B_{0, r}[g])
  & =\int_{B_{0, r}} \sqrt{\det g_{ij}} \tag 3.27\\
  & =\int_{B_{0, r}}\big[1-\frac 16 R_{ij}x_ix_j+o(|x|^2)\big]\\
  & =\frac 1n \omega_nr^n\big[1-\frac {R}{6(n+2)}r^2+o(r^2)\big],\\
  \endalign$$
where $R_{ij}$ and $R$ are respectively the Ricci curvature and
the scalar curvature in $g$. This is a contradiction when $R>0$ at
$0$, as the ratio $Q$ is a constant. Hence the scalar curvature of
$g$ vanishes almost everywhere. $\square$

\vskip10pt

\noo{\bf 3.4. End of proof of Theorem A}. From \S 3.3 and \S 3.4,
we see that if $(\M, g_0)$ is a compact manifold and there exists
a $k$-admissible function $w$ with singularity at some point $0$,
then $w$ has the asymptotic formula (3.17) and $w$ is smooth away
from $0$. The manifold $\M\back\{0\}$ equipped with the metric
$g=e^{-2w}g_0$ is a complete manifold with nonnegative Ricci
curvature, and satisfies furthermore the volume growth formula
$Q(r)\equiv 1$. Hence $(\M\back\{0\}, g)$ is isometric to the
Euclidean space [Cha]. Hence $(\M, g_0)$ is conformally equivalent
to the unit sphere $S^n$.

To finish the proof of Theorem A, it suffices to prove

\proclaim{Lemma 3.6} Let $(\M, g_0)$ be a compact manifold. If
$(\M, g_0)$ is not conformally equivalent to the unit sphere
$S^n$, then there exists $K>0$ such that if $w$ is a
$k$-admissible function,
$$\align
&\sup_\M w-\inf_\M w\le K, \tag 3.28\\
&|w(x)-w(y)|\le K|x-y|^{2-\frac nk}.\tag 3.29\\
\endalign$$
\endproclaim

\noo{\it Proof}. If (3.28) is not true, there exists a sequence of
$k$-admissible functions $w_m$ such that $\sup_\M w_m=0$ and
$\inf_\M w_m\to -\infty$. Suppose that $w_m(0)\to-\infty$. By the
H\"older continuity in \S3.1, we may assume that $e^{w_m}$
converges locally uniformly to $e^w$ in $\M\back\{0\}$. Obviously
$\lim_{x\to 0} w(x)=-\infty$. But from the above discussion, $(\M,
g_0)$ is conformally equivalent to the unit sphere $S^n$, which is
ruled out by out assumption. Hence (3.28) holds.

The H\"older continuity (3.29) follows from Lemma 3.1. $\square$

\vskip10pt

\noo{\bf 3.5. Remarks on the set $[g_0]_k$}. In this section we
prove some properties for $k$-admissible functions.

\proclaim{Lemma 3.7} If $w_1, w_2$ are smooth and $k$-admissible,
then $w=\max(w_1, w_2)$ is $k$-admissible.
\endproclaim

\noo{\it Proof}. It is convenient to consider the function
$u=e^w$. By approximation we suppose $u_1$ and $u_2$ are smooth
and $k$-admissible functions such that the eigenvalues $\lam(U)$
lie strictly in the open convex cone $\Ga_k$, where $U$ is the
matrix (1.16) with $u=u_1$ and $u_2$. Hence when $r>0$ is
sufficiently small, the eigenvalues of the matrix
$$U_r=\{u_{ij}-\frac {|\D u|^2}{2u_{x, r}}+uA_{g_0}\}\tag 3.30$$
lie in $\Ga_k$ for $u=u_1$ and $u_2$, where
$u_{x_0,r}=\inf_{B_{x_0, r}} u$.

Let $u=\max(u_1, u_2)$. Since $u_1, u_2$ are smooth function, $u$
is twice differentiable almost everywhere.  Let $\rho\in
C_0^\infty(\R^n)$ be a mollifier. In particular we choose $\rho$
to be a radial, smooth, nonnegative function, supported in the
unit ball $B_{0, 1}$, with $\int_{B_{0, 1}}\rho=1$.  Let
$$u_{[\eps]}(x)
 =\int_{B_{x, \eps}}\eps^{-n}\rho(\frac {|x-y|}\eps)
                        u(y) \sqrt{\det (g_0)_{ij}} dy\tag 3.31$$
be the mollification of $u$, where $B_{x, \eps}$ is the geodesic
ball. For each point $x$, using normal coordinates and the
exponential map, we have,  by (3.26),
$$\align
u_{[\eps]}(x)
 &=\int_{B_{0, 1}} \rho(y)u(x-\eps y)\sqrt{\det (g_0)_{ij}}\,dy \tag 3.32\\
 &=\int_{B_{0, 1}} \rho(y)u(x-\eps y)
    (1-\frac {\eps^2}{6}R_{ij}(x)y_iy_j+O(\eps^3))dy , \\
 \endalign $$
where $B_{0, 1}$ is the Euclidean space.  If $g_0$ is a flat
metric, we have
$$\align
\D u_{[\eps]} &= \int_{B_{0, 1}} \rho(y)\D u(x-\eps y)dy,\tag 3.33\\
\D^2 u_{[\eps]} &\ge \int_{B_{0, 1}} \rho(y)\D^2 u(x-\eps y)dy,\tag 3.34\\
 |\D u_{[\eps]}|^2
  &=[\int_{B_{0, 1}} \rho(y)\D u(x-\eps y)dy]^2\tag 3.35\\
  &\le \int_{B_{0, 1}} \rho(y)|\D u(x-\eps y)|^2dy.\\
  \endalign $$
Hence $u_{[\eps]}$ is $k$-admissible by (3.30). If $g_0$ is not
flat, by (3.32), an extra term of magnitude $O(\eps^2)$ arises.
Letting $\eps>0$ be sufficiently small and noting that the
eigenvalues of $U$ (with respect to $u_1$ and $u_2$) lie strictly
in the open set $\Ga_k$, we conclude again that $u_{[\eps]}$ is
$k$-admissible. $\square$

\proclaim{Corollary 3.8} Suppose $\phi$ is a smooth $k$-admissible
function on $\M$ with $\sigma_k(\lam(A_{g_\phi}))>f$, where
$g_\phi=e^{-2\phi} g_0\in [g_0]_k$ and $f$ is a smooth, positive
function. Let $w$ be the solution of
$$\align
\sigma_k(\lam(W))& =f\ \ \text{in}\ \ \Om, \tag 3.36\\
 w &=\phi\ \ \ \text{on}\ \ \pom,\\
 \endalign $$
where $W$ is given in (1.17), and $\Om$ is a smooth domain on
$\M$. Extend $w$ to $\M$ by letting $w=\phi$ on $\M-\Om$. Then $w$
is $k$-admissible.
\endproclaim

It was proved in [G] that (3.36) admits a solution $w$, smooth up
to the boundary. By the comparison principle we have $w>\phi$ in
$\Om$ and $\p_{\nu} (\phi-w)>0$ on $\pom$, where $\nu$ is the unit
outward normal. Hence we can extend $w$ to a neighbourhood of
$\Om$ such that it is $k$-admissible. Hence Corollary 3.8 follows
from Lemma 3.7.

\proclaim{Corollary 3.9} Consider the Dirichlet problem (3.36).
Suppose the set of sub-solutions $W_{sub}$ is not empty. Let
$$w(x)=\sup\{\phi(x)\ |\ \phi\in W_{sub}\}.\tag 3.37$$
If $w$ is bounded from above, then it is a solution to (3.36).
\endproclaim

By the interior a priori estimates [GW1, STW], the proof is
standard. Note that in Corollary 3.9, we allow $\Om$ to be the
whole manifold $\M$.

\vskip30pt

\centerline{\bf 4. Proof of Theorem C}

\vskip10pt

We divide the proof into three cases, according to $p<k$, $p=k$,
and $p>k$.

{\bf Case 1}: $p<k$. By (1.15), we can write equation (1.12) as
$$\sigma_k(\lam(W))=fe^{aw}, \tag 4.1$$
where
$$a=\frac 12(n-2)(k-p). \tag 4.2$$
For any given $k$-admissible function $w$, the functions $w+c$ and
$w-c$ are respectively a super and a sub solution of (4.1)
provided the constant $c$ is sufficiently large. By the a priori
estimates in [V2, GW1, STW] and the comparison principle, the
solution of (4.1) is uniformly bounded. When $a>0$, the linearized
equation of (4.1) is invertible. Hence by the continuity method,
there is a unique smooth solution to (4.1).

{\bf Case 2}: $p=k$. We prove that for any positive smooth
function $f$, there is a unique constant $\th>0$ such that the
equation
$$\sigma_k(\lam(W))=\th f \tag 4.3$$
has a solution. For $a>0$ small, let $w_a$ be the solution of
(4.1). Let $c_a=\inf w_a$. We write (4.1) in the form
$$\sigma_k(\lam(W_a))=(fe^{ac_a})e^{a(w_a-c_a)},\tag 4.4$$
where $W_a$ is the matrix (1.17) relative to $w_a$.  Assume
$g_0\in [g_0]_k$ so that $\lam(A_{g_0})\in\Ga_k$. Then at the
maximum point of $w_a$,
$$ \sigma_k(\lam(A_{g_0})\ge \sigma_k(\lam(W_a))\ge fe^{ac_a}.$$
At the minimum point of $w_a$,
$$ \sigma_k(\lam(A_{g_0})\le \sigma_k(\lam(W_a))= fe^{ac_a}.$$
Hence $e^{ac_a}$ is strictly positive and uniformly bounded as
$a\to 0$. By the a priori estimates [GW1, STW], where the
estimates depend only on $\inf (w_a-c_a)$,  we see that $w_a-c_a$
is uniformly bounded from above and sub-converges to a solution
$w_0$ of (4.3) with $\th=\lim_{a\to 0}e^{ac_a}$. By the maximum
principle it is easy to see that if $w'$ is another solution, then
necessarily $w'=w_0+const$; and furthermore (4.3) has no
($k$-admissible) solution for different $\th$.

{\bf Case 3}:  $p>k$. In this case we adopt the degree argument
from [W], see the proof of Theorem 5.1 there. Alternatively we can
also use the degree argument in \S 3 of [W]. We will study the
auxiliary problem
$$\sigma_k(\lam(V))=t(\delta_t+fv^p), \tag 4.6$$
where $t\ge 0$ is a parameter and $\delta_t$ is a positive
constant depending on $t$, $\delta_t=\delta_0\le 1$ when $t\le 1$
and $\delta_t=1$ when $t>2$, and $\delta_t$ is smooth and monotone
increasing when $1\le t\le 2$.

{\it Claim 1}. For any $t_0>0$, the solution of (4.6) is uniformly
bounded when $t\ge t_0$. Indeed, if there exists a sequence of
solutions $(t_j, v_j)$ of (4.6) such that $t_j\ge t_0$ and $\sup
v_j\to\infty$,  we have $m_j=\inf v_j\to\infty$ by (1.5). The
function $v'_j=v_j/m_j$ satisfies
$$\align
\sigma_k(\lam(V')) & \ge t_jf m_j^{p-k}(v'_j)^p)\\
    & \ge t_jf m_j^{p-k}\to\infty, \tag 4.7\\
    \endalign$$
where $V'$ is the matrix (1.11) relative to $v'$.  From (4.7) and
the comparison principle we have $\sup v_j'\to \infty$. Hence
$\inf v'_j\to\infty$ by (1.5), which contradicts to the definition
of $v_j'$.

Define the mapping $T_t$ so that for any $v_1\in C^2(\M)$,
$T_t(v_1)$ is the solution of
$$\sigma_k(\lam(V))=t(\delta_t+fv_1^p).\tag 4.8$$
Then a solution of (4.6) is a fixed point of $T_t$.

{\it Claim 2}. There is a solution of (4.6) when $t>0$ is small.
Indeed, for any smooth, positive function $\phi^*$, denote
$\Phi=\{\phi\in C^2(\M)\ |\ \phi<\phi^*\}$.  Then when $t>0$ is
small, $T(\Phi)$ is strictly contained in $\Phi$. Hence the degree
$\text{deg}(I-T_t, \Phi, 0)$ is well defined for $t\ge 0$ small.
Extend $T_t$ to $t=0$ by letting $T_t(v)=0$ for all $v$, so that
$T_t$ is also continuous at $t=0$. Hence
$$\text{deg}(I-T_t, \Phi, 0)=\text{deg}(I-T_0, \Phi, 0)=1. \tag 4.9$$
Hence $T_t$ has a fixed point in $\Phi$ for $t>0$ small.

{\it Claim 3}. Let $t^*=\sup\{t\ |\ (4.6) \text{ admits a
solution}\}$. Then $t^*$ is finite. Indeed, if $t^*=\infty$, there
is a sequence $t_j\to\infty$ such that (4.6) has a solution $v_j$.
We have obviously $m_j=\inf v_j\to\infty$, which is a
contradiction with Claim 1.

{\it Claim 4}. Equation (4.6) has a solution at $t=t^*$. Indeed,
let $t_j\nearrow t^*$ and $v_j$ be the corresponding solution of
(4.6). By claim 1, $v_j$ is uniformly bounded. Hence $v_j$
sub-converges to a solution $v^*$ of (4.6) with $t=t^*$.

Now we choose $\phi^*=v^*$ and define $\Phi$ as above. For any
$v_1\in\Phi$, let $v$ be the solution of (4.8). Since for any
$t\in (0, t^*)$, $v^*$ is a super-solution of (4.6). We have
$0<v<v^*$ by the maximum principle. Hence by (4.9),
$\text{deg}(I-T_t, \Phi, 0)=1$  for $t\in [0, t^*)$.

On the other hand, for any given $t_0>0$, since the solution of
(4.6) is uniformly bounded for $t\ge t_0$, the degree
$\text{deg}(I-T_t, B_R, 0)$ is well defined for $t\in (t_0,
t^*+1]$ for sufficiently large $R$, where $B_R=\{v\in C^2(\M)\ | \
v<R\}$. But when $t>t*$, (4.6) has no solution. Hence
$\text{deg}(I-T_t, B_R, 0)=0$. Hence for any $t\ge t_0$, (4.6) has
a solution $v\not\in\Phi$ with degree $-1$.

Let $v=v_{\delta_0}\not\in\Phi$ be a solution of (4.6) at $t=1$.
We have $\sup v>\inf v^*>0$. Let $\delta_0\to 0$. Since the
solution is uniformly bounded, it converges to a solution of
(1.12). This completes the proof. $\square$

From the above argument, we have the following extensions.

\proclaim{Theorem 4.1}  Let $(\M, g_0)$ be a compact $n$-manifold
not conformally equivalent to the unit sphere $S^n$. Suppose
$\frac n2<k\le n$ and $[g_0]_k\ne \emptyset$. Suppose there exists
a constant $c_0>0$ such that
$$\align
  \phi(x, t) & \ge c_0,\tag 4.10\\
  \lim_{t\to\infty} t^{-k} \phi(x, & t)  =\infty .\tag 4.11\\
 \endalign $$
Then there exists a constant $t^*>0$ such that the equation
$$\sigma_k(\lam(V))=t\phi(x, v) \tag 4.12$$
has at least two solutions for $0<t<t^*$, one solution at $t=t^*$,
and no solution for $t>t^*$.
\endproclaim

\proclaim{Theorem 4.2}  Let $(\M, g_0)$ be as in Theorem 4.1,
$\frac n2<k\le n$. Suppose $\phi>0$,
$$\lim_{t\to 0} t^{-k}\phi(x, t)  =0 ,\tag 4.13$$
and (4.11) holds. Then there exists a solution to (1.10).
\endproclaim

In the above theorems, we can also allow that the right hand side
depends on the gradient $\D v$. Furthermore, (4.11) and (4.13) can
be relaxed to
$$\align
&\lim_{t\to\infty} t^{-k} \phi(x,  t)  >\th, \tag 4.14\\
&\lim_{t\to 0}\, t^{-k}\phi(x, t)  <\th, \tag 4.15\\
\endalign $$
where $\th$ is the eigenvalue of (1.13) (with $f\equiv 1$). See
[W] for the Monge-Amp\'ere equation.

We remark that when $1\le k\le \frac n2$, Theorem C holds for
$p<k\frac {n+2}{n-2}$. Indeed, when $p\le k$, the proof of the
Cases 1 and 2 above also applies to the cases $1\le k\le \frac
n2$. When $k<p<k\frac {n+2}{n-2}$, by a blow-up argument and the
Liouville theorem [LL1], it is known that the set of solutions to
(4.6) is uniformly bounded. Hence by the above degree argument,
one also obtain the existence of solutions.

\proclaim{Theorem 4.3}  Let $(\M, g_0)$ be a compact $n$-manifold
with $[g_0]_k\ne \emptyset$, $1\le k\le n$. Then for any smooth,
positive function $f$ and any constant $p\ne k$, $p<k\frac
{n+2}{n-2}$, there exists a positive solution to the equation
(1.12). The solution is unique if $p<k$. When $p=k$, there exists
a unique constant $\th>0$ such that (1.13) has a solution. The
solution is unique up to a constant multiplication.
\endproclaim

Note that in Theorem 4.3 we allow that $(\M, g_0)$ is the unit
sphere.

\vskip20pt

\baselineskip=12.0pt
\parskip=2.0pt

\Refs\widestnumber\key{ABC}

\item {[A1]} T. Aubin, Equations diff\'erentielles non lin\'eaires
       et probl\`me de Yamabe concernant la courbure scalaire.
       J. Math. Pures Appl. (9) 55 (1976), 269--296.

\item {[A2]} T. Aubin,
       Some nonlinear problems in Riemannian geometry,
       Springer, 1998.

\item {[CNS]} L.A. Caffarelli, L. Nirenberg, and J. Spruck,
       Dirichlet problem for nonlinear second order
       elliptic equations III.
       Functions of the eigenvalues of the Hessian,
       Acta Math. 155(1985), 261--301.

\item {[C]} J.-G. Cao,
       The existence of generalized isothermal coordinates
       for higher-dimensional Riemannian manifolds,
       Trans. Amer. Math. Soc. 324(1991), 901--920.

\item {[CGY1]} A. Chang, M. Gursky, P. Yang,
       An equation of Monge-Am\`ere type in conformal geometry,
       and four-manifolds of positive Ricci curvature,
       Ann. of Math. (2) 155(2002), 709--787.

\item {[CGY2]} A. Chang, M. Gursky, P. Yang,
       An a priori estimate for a fully nonlinear equation on
       four-manifolds,
       J. Anal. Math. 87 (2002), 151--186.

\item {[CHY]} Alice Chang, Z.-C. Han, P. Yang,
       Classification of singular radial solutions to the
       $\sigma_k$-Yamabe equation on annular domains, preprint.

\item {[Cha]} I. Chavel,
       Riemannian geometry---a modern introduction,
       Cambridge Univ. Press, 1993.

\item {[Ch]} K.S. Chou (K. Tso),
      On a real Monge-Ampere functional,
      Invent. Math. 101(1990), 425--448.

\item {[CW]} K.S. Chou and X-J. Wang,
       A variational theory of the Hessian equation,
       Comm. Pure Appl. Math. 54 (2001), 1029--1064.

\item {[GeW]} Y. Ge and G. Wang,
       On a fully nonlinear Yamabe problem, preprint.

\item {[GT]}  D. Gilbarg and N.S. Trudinger,
       Elliptic partial differential equations of second order,
       Springer, 1983.

\item {[G]} B. Guan,
       Conformal metrics with prescribed curvature functions on
       manifolds with boundary, \newline preprint.

\item {[GW1]} P. Guan and G. Wang,
       Local estimates for a class of fully nonlinear equations
       arising from conformal geometry,
       Int. Math. Res. Not. (2003), 1413--1432.

\item {[GW2]} P. Guan and G. Wang,
        A fully nonlinear conformal flow on locally conformally
        flat manifolds,
        J. Reine Angew. Math. 557 (2003), 219--238.

\item {[GVW]} P. Guan, J. Viaclovsky, and G. Wang,
      Some properties of the Schouten tensor and applications
      to conformal geometry,   Trans. Amer. Math. Soc., 355(2003),
      925-933.

\item {[Gu]} M. G\"unther,
       Conformal normal coordinates,
        Ann. Global Anal. Geom., 11(1993), 173--184.

\item {[GV1]} M. Gursky and J.Viaclovsky,
       Prescribing symmetric functions of the eigenvalues of the
       Ricci tensor, arXiv:math.DG/0409187.

\item {[GV2]} M. Gursky and J.Viaclovsky,
       Convexity and singularities of curvature equations in
       conformal geometry, arXiv:math.DG/0504066.

\item {[LP]} J.M. Lee and T.H. Parker,
        The Yamabe problem, Bull. Amer. Math. Soc. 17(1987), 37--91.

\item {[LL1]} A. Li and Y.Y. Li,
       On some conformally invariant fully nonlinear equations.
       Comm. Pure Appl. Math. 56 (2003), 1416--1464.

\item {[LL2]} A. Li and Y.Y. Li,
       On some conformally invariant fully nonlinear equations
       I\!I, Liouville, Harnack, and Yamabe, preprint.

\item {[STW]} W.M. Sheng, N.S. Trudinger, X.-J. Wang,
       The Yamabe problem for higher order curvatures,
       preprint.

\item {[S]} R. Schoen,
       Conformal deformation of a Riemannian metric to constant
       scalar curvature,
       J. Diff. Geom. 20(1984), 479--495.

\item {[SY]} R. Schoen and S.T. Yau,
       Lectures on Differential geometry.
       International Press, 1994.

\item {[T]} N.S. Trudinger,
       On Harnack type inequalities and their application
       to quasilinear elliptic equations,
       Comm. Pure Appl. Math., 20(1967), 721-747.

\item {[TW1]} N.S. Trudinger and X-J. Wang,
        Hessian measures I,
        Topol. Methods Nonlinear Anal. 10 (1997), 225--239.

\item {[TW2]} N.S. Trudinger and X-J. Wang,
       Hessian measures I\!I,
       Ann. of Math. (2) 150 (1999), 579--604.

\item {[V1]} J. Viaclovsky, Conformal geometry, contact geometry,
      and the calculus of variations.
      Duke Math. J. 101 (2000), no. 2, 283--316.

\item {[V2]} J. Viaclovsky,
       Estimates and existence results for some fully nonlinear
       elliptic equations on Riemannian manifolds,
       Comm. Anal. Geom. 10 (2002), 815--846.

\item {[W]} X.-J. Wang,
      Existence of multiple solutions to the equations of
      Monge-Amp\'ere type,
      J. Diff. Eqns, 100(1992), 95-118.

\endRefs

\enddocument
\end